\documentclass{birkjour}

\usepackage{
amsfonts,
latexsym,
amssymb,
amsmath,
amsthm,
enumerate,
verbatim,
mathrsfs,
epigraph}

\usepackage{url}

\newcommand{\labbel}{\label}

\newcommand{\uu}{{\hspace {1 pt} =  \hspace {1 pt}}}

\newtheorem{theorem}{Theorem}

\newtheorem{corollary}[theorem]{Corollary}

\newtheorem*{theorem*}{Theorem}
\newtheorem*{corollary*}{Corollary}

\theoremstyle{definition}
\newtheorem{definition}[theorem]{Definition}

\theoremstyle{remark}
\newtheorem{remark}[theorem]{Remark}

\newcommand{\brfrt}{\hspace{0 pt}}

 \allowdisplaybreaks[4]

\begin{document}

\title[Gumm and alvin levels]
{The Gumm level equals the alvin level in congruence distributive varieties}

\author{Paolo Lipparini} 
\address{Dipartimento di Matematica\\Viale della  Ricerca
 Scientifica\\Universit\`a di Roma ``````Tor Vergata'''''' 
\\I-00133 ROME ITALY}

\email{lipparin@axp.mat.uniroma2.it}

\subjclass{Primary 08B10}

\keywords{ Gumm terms,  alvin terms, J{\'o}nsson terms,
congruence distributive variety, congruence modular
variety, congruence identity}

\date{\today}

\thanks{\today
\\
Work performed under the auspices of G.N.S.A.G.A. Work 
partially supported by PRIN 2012 ``Logica, Modelli e Insiemi''.
The author acknowledges the MIUR Survival Department Project awarded to the
Department of Mathematics, University of Rome Tor Vergata, CUP
E83C18000100006.
}

\begin{abstract}
Congruence modular and congruence distributive 
varieties are  characterized by the existence of sequences of 
Gumm and J{\'o}nsson  terms, respectively.
Such sequences  
have variable lengths, in general.

It is immediate from the above paragraph that
 there is a variety with  Gumm terms but without
 J{\'o}nsson  terms. 
We prove the quite unexpected result 
that, on the other hand, if some variety has both kinds of terms,
then the minimal lengths of the sequences differ at most by $1$.

It follows that every $r$-modular  congruence distributive 
 variety is 
$r^2-r+2$-distributive. 
\end{abstract}

\maketitle

An \emph{algebra} (short for \emph{algebraic system})
is a nonempty set endowed with a family of operations.
A \emph{variety} is a class of algebras
of the same type which is  definable by a set of equations.
A \emph{congruence} on some algebra is the kernel of some
homomorphism, equivalently, a compatible equivalence relation.
 The set of congruences on some algebra has a lattice structure.
An algebra is \emph{congruence modular} 
if its lattice of congruences is modular and 
a variety is congruence modular 
if so are all of its members.
Congruence distributivity is defined in a similar manner. 
Congruence modular  varieties 
include the varieties of groups, of rings, of quasigroups, as well as
all congruence distributive varieties.
Congruence distributive varieties 
include the varieties of lattices and of Boolean algebras.

Of course, there have been 
interactions between the
general theory of algebras and more specific 
kinds of algebraic systems \cite{M}.
Quite unexpectedly, interesting connections 
recently emerged with the theory of 
computational complexity, in particular, the
algebraic approach to the Constraint Satisfaction Problem. 
In a nutshell,
well-behaved classes 
of algebras first studied only for their algebraic properties---and
the identities such classes satisfy---turned out to
 correspond to
algorithmic results for CSP
over constraint languages.
 See, e.~g.,
\cite{BKW,JKN} for a survey. 
Congruence distributivity and 
congruence modularity 
are among the first studied and most important
conditions providing `good algebraic
 structure'. 
We show that two well-known
and intensively studied characterizations 
of the above conditions are 
 equivalent even as far as length is concerned,
 provided they are both applicable.

In more detail, a \emph{term} is, roughly, a word obtained by 
composition from the basic operations of a variety. 
Both congruence distributivity and congruence modularity
are characterized by the existence
of  finite sequences of terms satisfying appropriate equations.
In general, the lengths of  such sequences
are variable. Characterizations of this kind are called Maltsev 
conditions. 
 See Theorems \ref{alvin}, \ref{gumm} 
and Definition \ref{gummdef}  below for specific examples.
Many problems are still open about the relationship
among the relative lengths of different sequences
characterizing different conditions.
Such problems date back at least to
\cite[p.\ 173]{D}. See, e.~g., 
\cite{CHL,FV,CV,adjt,LTT,jds,ntcm,dsh,T} 
for other problems, comments and
results. The reader can find  further references in the quoted works.
See also  Remark \ref{gummday} below.

Gumm terms and alvin terms, defined below, characterize congruence modularity
and  congruence distributivity, respectively.
In particular, there is a variety with a sequence of Gumm terms but without
a sequence of alvin terms.
In this note we prove the quite unexpected result 
that, within a variety and \emph{as soon as we have both kinds of terms}, 
the minimal lengths of the sequences are identical.
While it is well established that Gumm
terms are relevant for the study of congruence modular 
varieties \cite{FMK,G,L,LTT,T},
it is quite surprising to discover
that they have a deep direct influence on 
congruence distributive varieties.
Of course, a congruence distributive variety is also congruence modular,
but it should be expected that in a congruence distributive 
variety ``all the work'' 
is done by an alvin sequence.
On the contrary, we show that the weaker notion of a Gumm sequence
has applications, if the sequence is  sufficiently short.

Congruence modularity 
can be characterized by a different set of terms
formerly  introduced by A.\ Day.
The problem of the relationship between the 
minimal lengths of Day and Gumm 
sequences in a congruence modular variety 
is not 
completely solved yet.
The results presented in this note stress
the importance of the above problem.
Anyway, we get  the remarkable corollary  that
an $r$-modular congruence distributive variety is
$r^2-r+2$-distributive, though we do not 
know how far this bound can be improved. 
In passing, we get new proofs, usually simpler, usually 
providing better bounds, of results 
from \cite{CHL,jds,ntcm,T}. 

Let us now recall the main definitions together with
some needed classical results.
See, e.~g., \cite{FMK,G,CV,dsh,MMT}  
for further undefined notions and full formal details.

\begin{theorem} \labbel{alvin}
\cite{JD,MMT} 
A variety $\mathcal V$ is congruence distributive 
if and only if  there is some natural number $n$ 
for which one of the following equivalent 
conditions holds.
\begin{enumerate}[(1)]  
 \item
$\mathcal V$ has a sequence $t_0, \dots, t_n$ of \emph{alvin terms}, 
that is,  terms such that 
 the following equations are satisfied in
all algebras in $\mathcal V$.
\begin{align}
\labbel{a1}    \tag{A1}
  x  &=t_h(x,y, x), 
&&\text{for } 0 < h < n,   
\\
\labbel{a2}    \tag{A2}     
  x&= t_0(x,y,z),  
 \\
\labbel{a3}    \tag{A3}     
t_{h}(x,z,z) &=
t_{h+1}(x,z,z), 
&&\text{for $h$ even, } 0 \leq h < n,   
\\ 
\labbel{a4}    \tag{A4}     
 t_{h}(x,x,z) &=
t_{h+1}(x,x,z),
&&\text{for $h$ odd, } 0 \leq h <  n,   
\\
\labbel{a5}    \tag{A5}     
 t_{n}(x,y,z)&=z. 
\end{align}   
\item
The following inclusion holds
in every algebra $\mathbf A \in \mathcal V$ 
for all congruences $\alpha$, $\beta$ and $\gamma$ 
of $\mathbf A$. 
\begin{equation}
\labbel{cd}    \tag{CD}      
\alpha ( \beta  \circ  \gamma ) 
\subseteq 
\alpha \gamma   \circ 
 \alpha \beta \circ  
 {\stackrel{n}{\dots}} 
 \end{equation}     
  \end{enumerate} 
 \end{theorem} 

The notation we use in congruence identities
like \eqref{cd} above goes as follows.
Juxtaposition denotes intersection, in particular, meet of 
congruences. 
Join in congruence lattices is denoted by $+$.
 Composition
of binary relations 
is denoted by 
 $\circ$ 
and  
$R  \circ S \circ  {\stackrel{n}{\dots}}  $ 
 denotes 
$R  \circ S \circ R \circ S \dots $
with $n$ factors, that is, $n-1$ occurrences of 
$\circ$. In the above 
notation factors of the form $\alpha \beta $ are always counted as
\emph{one} factor.
It is formally convenient to allow 
the extreme cases $n=0$
and $n=1$.
We set    
$R  \circ S \circ  {\stackrel{1}{\dots}} =R $
and 
$R  \circ S \circ  {\stackrel{0}{\dots}}  $
to be the minimal congruence of the algebra under consideration.
We also let 
$ R ^{m} = R  \circ R \circ  {\stackrel{m}{\dots}}  $
Exponentiation
ties more than 
any other operator;
juxtaposition comes next in tying force.

The equivalence of (1) and (2)
in Theorem \ref{alvin} is implicit in 
\cite{JD} and is an immediate application 
of the algorithm described in 
\cite{Pal,W}. See also 
\cite{LTT,jds,T} for further comments and related results. 

The classical characterization \cite{JD} of congruence distributive
varieties 
involves \emph{J{\'o}nsson terms}, which are defined as in 
clause (1) above with even and odd exchanged in 
\eqref{a3} - \eqref{a4}.  Of course, the definitions 
are completely equivalent if we are not concerned 
with the exact value of $n$. 
However, for a fixed even $n$, the alvin and the J{\'o}nsson 
conditions are not equivalent \cite{FV}, though it is obvious that
the minimal lengths of the sequences differ at most by one. 
As mentioned, we shall show that 
the minimal lengths of an alvin and a Gumm
sequence coincide in a congruence distributive variety.
Since the minimal length of a J{\'o}nsson 
sequence might differ by $1$, 
 we get an
exact correspondence 
only if we deal with the alvin condition.

\begin{definition} \labbel{gummdef}   
A sequence
of \emph{Gumm terms}
for a variety is a sequence $t_0, \dots, t_n$,
for some $n$, of terms 
satisfying
the equations \eqref{a2} - \eqref{a5}  
in Clause (1) in Theorem \ref{alvin} above,
as well as
\begin{align}
\labbel{g1}    \tag{G1}
  x  &=t_h(x,y, x), 
&&\text{for }  1 < h < n.   
\end{align}
In other words, 
Gumm terms satisfy 
the equations \eqref{a1} - \eqref{a5}, 
except possibly for the equation 
$x=t_1 (x,y,x)$. 
 \end{definition}

Notice that 
$x=t_1(x,x,x)$ still holds in the case of Gumm terms, 
by \eqref{a2} and \eqref{a3}. 

\begin{theorem} \labbel{gumm}
\cite{G1,G}  A variety $\mathcal V$ is congruence modular
if and only if 
$\mathcal V$ has a sequence of Gumm terms.
 \end{theorem}

A variety is \emph{$n$-alvin} (\emph{$n$-Gumm, $n$-distributive})
 if it has a sequence $t_0, \dots, \allowbreak t_n$ 
of alvin (Gumm, J{\'o}nsson) terms.
The \emph{alvin (Gumm) level} 
of a congruence distributive (modular) variety $\mathcal V$  
is the minimal $n$
such that $\mathcal V$ is $n$-alvin ($n$-Gumm).
Our main result about the above levels is the following theorem.
The proof shall be given after some auxiliary results
of independent interest.

\begin{theorem} \labbel{=}
If $\mathcal V$ is a congruence distributive variety,
then the alvin level of $\mathcal V$ 
is equal to the Gumm level of $\mathcal V$.
 \end{theorem}

Notice that the sequence 
$t_0, \dots, t_n$
actually contains $n+1$ terms;
moreover, 
the two ``outer''  terms 
$t_0 $ and $ t_n$ 
are  projections, hence the number of nontrivial
terms is $n-1$. Thus any definition
of the levels has a somewhat conventional 
nature.
Here we follow by analogy the classical 
and universally adopted convention 
concerning $n$(-J{\'o}nsson)-distributivity.

If either $n=0$
or $n=1$ in
Theorem \ref{alvin}
or in Definition \ref{gummdef},
then we get a condition which is satisfied only
by trivial varieties 
with just one-element algebras.
Indeed, say, for $n=1$
we get 
$x= t_0(x,y,y) = t_1(x,y,y) =y$
from \eqref{a2}, \eqref{a3} and \eqref{a5}.  
However, it is formally convenient 
to consider the above trivial cases, too, 
otherwise the alvin and Gumm levels would
be undefined in the case of trivial varieties.
Of course, the levels trivially coincide
(and both equal $0$) in such trivial cases.

According to Definition \ref{gummdef},  
Gumm terms can be seen as a ``defective''
variant of
alvin terms. 
 While this idea is sometimes useful,
Gumm terms have a more important 
and fruitful interpretation as terms which ``compose''  
the classical Maltsev conditions 
for congruence permutability and distributivity.
See \cite{G1,G,ntcm,T}  for a more detailed discussion.
In some respects, 
the main point of the present note is to 
stretch 
 the ``permutable''
side of the term $t_1$ to the extreme limit.  
Notice also that the present definition 
is slightly different from the original definition
by H.-P. Gumm from \cite{G1,G}.
The present definition allows for a finer counting of the 
number of terms. See \cite[p.\ 12]{jds}.
To the best of our knowledge, the present definition 
first appeared in \cite{LTT,T}.  
 Notice that the indexing of terms, whatever the definition,
is different from the present one in most of the quoted papers,
including works by the present author.
 The indexing here is  intended to stress the similarity
between Gumm and alvin terms.

We now present our main new applications
of Gumm terms. The next theorem  works
for any congruence modular variety.
A \emph{tolerance}
on some algebra $\mathbf A$ 
is a reflexive and symmetric binary
 compatible relation on $A$.
In other words, a  tolerance is like a congruence, except that 
transitivity is not required.
Recall the notational conventions
established shortly after
the statement of Theorem \ref{alvin}. 
For  convenience, we shall 
frequently write $a \mathrel R b$
in place of $(a,b) \in R$ and we shall also concatenate
the above notation, e.~g., 
 $a \mathrel R b \mathrel { S} c $
means both $(a,b) \in R$ and $(b, c) \in S$.

\begin{theorem} \labbel{reduction}
Suppose that $n \geq 1$
and $\mathcal V$ is a  variety 
with a sequence 
$t_0, \dots, t_n$ of Gumm (alvin) terms.
Then, for every $m \geq 1$,  $\mathcal V$ has a sequence 
$s_0, \dots, s_n$ of Gumm (alvin) terms
such that $s_1$ satisfies the 
following additional property.

  \begin{enumerate}    \item []
  \begin{enumerate}    \item    
[(T$_m$)] For every $\mathbf A \in \mathcal V$ 
and every tolerance $\Theta$ on $\mathbf A$,
if $a,c \in A$ and $ a \mathrel { \Theta ^m} c$, 
then $ a \mathrel { \Theta } s_1 (a,a,c) $. 
 \end{enumerate} 
\end{enumerate} 
\end{theorem}

\begin{proof}
The property 
(T$_1$)
is trivially satisfied by the 
term $t_1$ of the  
original sequence,
since if  
$ a \mathrel { \Theta } c$, 
then $ a = t_1 (a,a,a) \mathrel { \Theta } t_1 (a,a,c) $,
by equations \eqref{a2} - \eqref{a3}
and since $\Theta$ is reflexive and compatible.

We shall prove that if $m \geq 1$
and  $\mathcal V$ has a sequence 
$s_0, \dots, s_n$ of Gumm (alvin) terms such that 
$s_1$ satisfies  
(T$_m$), then 
$\mathcal V$ has a sequence 
$s^*_0, \dots, s^*_n$ of Gumm (alvin) terms such that 
$s^*_1$ satisfies  
(T$_{m+1}$). The theorem follows by induction on $m$. 
Define 
\begin{equation}\labbel{*}  \tag{*} 
s_h^*(x,y,z) = s_h (x, s_h (x,y,y), s_h (x,y,z)), \quad 
\text{ for } h=0, \dots, n. 
   \end{equation}     

If
$s_0, \dots, s_n$ is  a sequence of Gumm terms,
then 
\begin{align*}\labbel{dim}
  s^*_h(x,y, x) &= s_h (x, s_h (x,y,y), s_h (x,y,x)) =
s_h (x, s_h (x,y,y), x) = x \ (h >1),   
\\
  s^*_0(x,y, z) &= s_0 (x, s_0 (x,y,y), s_0 (x,y,z)) = x,  
 \\
 s^*_h(x,z,z) &= s_h (x, s_h (x,z,z), s_h (x,z,z)) =
\\
& = s_{h+1} (x, s_{h+1} (x,z,z), s_{h+1} (x,z,z)) =
s^*_{h+1} (x,z,z)
\text{  ($h$ even),}   
\\
  s^*_h(x,x,z) &= s_h (x, s_h (x,x,x), s_h (x,x,z))
= s_h (x, x, s_{h+1} (x,x,z)) =
\\
&=s_{h+1}(x, x, s_{h+1} (x,x,z)) = s^*_{h+1} (x,x,z)
\text{ ($h$ odd),}
\\
  s^*_n(x,y, z) &= s_n (x, s_n (x,y,y), s_n (x,y,z)) = 
s_n (x,y,z) = z,
\end{align*}   
thus
$s^*_0, \dots, s^*_n$ is  a sequence of Gumm terms, as well.
The case of alvin terms is identical, just let $h$
be arbitrary in the first displayed line. 

Suppose now that 
$\mathbf A $ belongs to $  \mathcal V$, 
$\Theta$ is a tolerance on $\mathbf A$,
 $a,c \in A$ and $ a \mathrel { \Theta ^{m+1}} c$.
Thus there is $b \in A$
such that  $ a \mathrel { \Theta ^{m}} b \mathrel  \Theta c$.
Then 
\begin{equation*}
    a= s_1 (a, \underline{s_1(a,a,b)}, s_1(a,a, \underline{b})) \mathrel { \Theta } 
s_1 (a, \underline{a}, s_1(a,a, \underline{c})) =  s_1^*(a,a,c) 
 \end{equation*}
 by
\eqref{a2} - \eqref{a3}, 
where we have underlined the 
$\Theta$-related elements. 
We have used the
assumptions that
 $\Theta$ is reflexive and compatible
and that $s_1$ satisfies 
(T$_m$), thus $s_1(a,a,c) \mathrel \Theta a$,
since $\Theta$ is symmetric.  
 \end{proof}  

Compare  \eqref{*} with 
the definition shortly before 
\cite[Observation 10.1 on p.\ 64]{G}.  
Compare also \cite[Section 3]{adjt}. 
The position
\eqref{*} 
sends (directed, Pixley  \cite{adjt}) J{\'o}nsson terms to 
(directed, Pixley)  J{\'o}ns\-son 
terms,  as well.
In particular, Theorem \ref{reduction}
applies to
 \emph{Pixley terms} in the sense of \cite{adjt}.
Theorem \ref{reduction} applies also to
 \emph{directed Gumm terms} \cite{adjt},
provided either we define directed Gumm terms 
in a specular way,
or else we consider the variant of
  (T$_m$) whose conclusion
asks for $  s_{n-1} (a,c,c) \mathrel { \Theta } c$.

The complexity of the terms constructed in the 
 proof of Theorem \ref{reduction} is not optimal. 
Moreover, using some additional arguments,
it is possible to prove 
the version of Theorem \ref{reduction} 
in which the tolerance $\Theta$ 
is replaced by a reflexive and compatible relation.
See the appendix.
We shall not need the above generalizations here,
hence  we have favored ease 
over generality.

 Only part (1) of the following corollary
 shall be used  in order to prove 
Theorems \ref{=} and \ref{diLTT}.

\begin{corollary} \labbel{coru}
Suppose that $\ell \geq 1$,
$n \geq 2$ and  $\mathbf A$ is an algebra
belonging to a variety  with Gumm terms
$t_0, \dots, t_n$. If
$\alpha$, $\beta$, $\gamma$ are congruences and
 $\Theta$, $\Psi$ are tolerances on $\mathbf A$, 
then the following inclusions hold.

  \begin{enumerate}[(1)]  
  \item 
$( \beta \circ \gamma )(\alpha \beta + \alpha \gamma )
 \subseteq
\alpha \gamma \circ \alpha \beta 
\circ  {\stackrel{n}{\dots}}$ (here possibly $n=0$ or $n=1$),
  \item 
$ \alpha ( \beta \circ \gamma \circ  {\stackrel{\ell}{\dots}}    )
 \subseteq
 \alpha ( \beta \circ \gamma)( \gamma \circ \beta ) \circ
(\alpha \beta  \circ \alpha \gamma  
\circ  {\stackrel{k}{\dots}})$, for $k= (n-2)( \ell -1) +1$, 
\item
$( \beta \circ \gamma  \circ  {\stackrel{\ell}{\dots}})
(\alpha \beta + \alpha \gamma )
 \subseteq
\alpha \gamma \circ \alpha \beta 
\circ  {\stackrel{k}{\dots}} \,$, for $k= (n-2)( \ell -1) + 2$, 
\item
$\Psi \Theta^ \ell \subseteq 
(\Psi \Theta)^{\ell(n-2)+1} $.
  \end{enumerate}
 \end{corollary}

 \begin{proof} 
(1) As we mentioned before, if 
$n \leq 1$, then we are in a trivial variety,
hence the conclusion holds. 
So let us suppose $n \geq 2$.

Let $\Theta$ 
be the  tolerance 
$ (\alpha \beta \circ \alpha \gamma)
(\alpha \gamma  \circ \alpha \beta )$
so that $\Theta$ contains
both $\alpha \beta $ and $\alpha \gamma $.
Suppose that $ (a,c) \in ( \beta \circ \gamma )(\alpha \beta + \alpha \gamma )$.
Using the classical characterization of the join of two
congruences, from $(a,c)\in  \alpha \beta + \alpha \gamma $
we get that there is some $m$ (depending on $a$ and  $c$, in 
general) such that 
$(a,c)\in  \alpha \beta \circ  \alpha \gamma \circ  {\stackrel{m}{\dots}} 
\subseteq \Theta^m    $.
By Theorem \ref{reduction},
we have Gumm terms 
$s_0, \dots, s_n$ such that 
$s_1$ satisfies 
(T$_{m}$),
hence 
$(a, s_1(a,a,c)) \in \Theta \subseteq \alpha \gamma \circ \alpha \beta  $,
by \eqref{a4} and  \eqref{a5}.
If $n=2$, then 
$s_1(a,a,c) = c$ and we are done.  
If $n >2$, it follows from a classical argument in 
\cite{JD} 
(or see the proof of  (2) below)
that    
$ (s_2(a,a,c), c) \in \alpha \beta \circ \alpha  \gamma 
\circ  {\stackrel{n-1}{\dots}}\, $,
since 
$(a,c) \in \alpha ( \beta \circ \gamma )$.
By
\eqref{a4} we get
$s_1(a,a,c) = s_2(a,a,c)$, hence
$(a,c) \in (\alpha \gamma \circ \alpha \beta)
\circ (\alpha \beta \circ \alpha  \gamma 
\circ  {\stackrel{n-1}{\dots}})=
\alpha \gamma \circ \alpha \beta 
\circ  {\stackrel{n}{\dots}}\,$,
since composition
of relations is associative and
$\alpha \beta $
is a congruence, hence transitive. 
 
(2) Let $\Psi$ 
be the  tolerance 
$ ( \beta \circ  \gamma)
(\gamma  \circ  \beta )$, 
thus $\Psi$ contains
both $ \beta $ and $ \gamma $.
Suppose that $ (a,c) \in 
 \alpha ( \beta \circ \gamma \circ  {\stackrel{\ell}{\dots}} ) $,
hence 
$ (a,c) \in  
\Psi ^ \ell $, 
thus  
$(a, s_1(a,a,c)) \in \Psi $,
for the appropriate term $s_1$ given by 
Theorem \ref{reduction}.
We have
$a = s_1(a,a,a) \mathrel { \alpha  }  
 s_1(a,a,c)$,
so that 
$(a, s_1(a,a,c)) \in \alpha ( \beta \circ  \gamma)
(\gamma  \circ  \beta )$.
If $n=2$, we are done, actually, we have $k=0$
and we can save one factor.

Classical  arguments show that
$ (s_2(a,a,c),c) \in \alpha \beta   \circ \alpha \gamma   
\circ  {\stackrel{k}{\dots}}  \, $,
for $k= (n-2)( \ell -1) +1$. 
Say, if $\ell$  is even 
and $a \mathrel { \beta } b_1 \mathrel { \gamma } b_2 
\mathrel { \beta } b_3 \mathrel { \gamma  } \dots 
\mathrel { \beta } b _{\ell -1} \mathrel \gamma c $,
then 
\begin{multline*} 
s_2(a,a,c) \mathrel { \beta } s_2(a,b_1,c) \mathrel { \gamma } s_2(a,b_2,c)
\mathrel { \beta } s_2(a,b_3,c) \mathrel { \gamma }
  \dots
\\
\dots \mathrel { \beta } s_2(a,b _{\ell -1},c) \mathrel { \gamma } 
s_2(a,c,c) = s_3(a,c,c) \mathrel { \gamma  } 
s_3(a,b _{\ell -1},c) \mathrel  \beta  \dots  
\end{multline*}   
Notice that then 
$s_2(a,b _{\ell -1},c) \mathrel { \gamma  } s_3(a,b _{\ell -1},c)$,
since $\gamma$ is transitive. 

All the  elements 
in the above chain are $\alpha$-connected by 
\eqref{g1}, e.~g.,
$s_2(a,a,c) \mathrel \alpha  s_2(a,a,a) =
a= s_2(a,b_1,a) \mathrel \alpha s_2(a,b_1,c)  $, etc.
Notice that we do not need the identity
$x=s_1(x,y,x)$ in order to show
 $a = s_1(a,a,a) \mathrel \alpha  s_1(a,a,c) = s_2(a,a,c) $,
the equations \eqref{a2} - \eqref{a4} are enough. 
The conclusion of (2) follows from
$s_1(a,a,c) = s_2(a,a,c)$, as in (1).
  
The proof of  (3) merges the arguments 
in (1) and (2).
As in the proof of (1)
and for the same definition of $\Theta$,
we get
$(a, s_1(a,a,c)) \in \Theta \subseteq \alpha \gamma \circ \alpha \beta  $,
for the appropriate $s_1$ 
given by Theorem \ref{reduction}.
The rest goes as in the last part of the proof of (2).
Finally, notice that here, as in (1), two adjacent occurrences 
of $\alpha \beta $ join into one, by transitivity. 

(4) As in the above arguments, 
if $(a,c) \in \Psi \Theta^ \ell $, 
we have 
$a \mathrel { \Theta }  s_1(a,a,c)$,
for some appropriate $s_1$, and 
 $a =   s_1(a,a,a) \mathrel { \Psi }  s_1(a,a,c)$,
thus 
$a \mathrel { \Psi \Theta }  s_1(a,a,c)$.
If $\Psi$ is a congruence, the  arguments in (2)
 give  
$( s_2(a,a,c), c) \in (\Psi \Theta)^{\ell(n-2)}$
(here we are not allowed to use transitivity in order to get
a better value).
The conclusion follows again from
$s_1(a,a,c) = s_2(a,a,c)$. 
In order to deal with the case when 
$\Psi$ is a tolerance, it is enough to use an additional argument from
\cite{CH}. E.~g.,
\begin{equation*} 
s_2(a{,\hspace{0.2 pt}}a{,\hspace{0.2 pt}}c) \uu 
s_2(s_2(a{,\hspace{0.2 pt}}a{,\hspace{0.2 pt}} \underline{c}) 
{,\hspace{0.2 pt}}b_1{,\hspace{0.2 pt}}
s_2( \underline{a}{,\hspace{0.2 pt}}a{,\hspace{0.2 pt}}c) )
{\hspace{1 pt}\Psi\hspace{1 pt}}
  s_2(s_2(a{,\hspace{0.2 pt}}a{,\hspace{0.2 pt}} \underline{a}) 
{,\hspace{0.2 pt}}b_1{,\hspace{0.2 pt}}
s_2( \underline{c}{,\hspace{0.2 pt}}a{,\hspace{0.2 pt}}c) )
\uu s_2(a{,\hspace{0.2 pt}}b_1{,\hspace{0.2 pt}}c),
  \end{equation*}    
by \eqref{g1}. See \cite{CH}
or the proof of \cite[Proposition 3.1]{jds}
for full details. 
\end{proof}  

Letting $\ell$  vary in  (2) above
and using \cite{G1},  we get another proof of the result
from \cite{T} that a variety $\mathcal V$ is congruence modular 
if and only if the congruence identity
$ \alpha ( \beta + \gamma ) = \alpha ( \beta \circ \gamma )
\circ ( \alpha \beta + \alpha \gamma )$ holds in $\mathcal V$.
Clause (2) generally gives the best known bound for 
$\alpha ( \beta \circ \gamma \circ  {\stackrel{\ell}{\dots}}    )$, 
to date.
For example, it follows from
\cite[Corollary 2.2]{jds}
that an $m+1$-distributive variety
satisfies
$\alpha ( \beta \circ \gamma \circ  {\stackrel{\ell+1}{\dots}} )
\subseteq  
\alpha \beta \circ \alpha \gamma  \circ  {\stackrel{m \ell +1}{\dots \dots}} $    
In case $m$ is even (thus $m+1$ is odd) 
the property of being $m+1$-distributive is 
equivalent to $m+1$-alvin, hence 
Corollary \ref{coru}(3) provides the improved 
inclusion
$\alpha ( \beta \circ \gamma \circ  {\stackrel{\ell+1}{\dots}} )
\subseteq  
\alpha \gamma  \circ \alpha \beta \circ 
 {\stackrel{(m-1) \ell +2}{\dots \dots}} $

Moreover, using (4), we  get bounds
for expressions of the form
\begin{equation*}
   ( \beta  _{11}  \circ \beta _{12} \circ \beta _{13} \dots  )
( \beta  _{21}  \circ \beta _{22} \circ \beta _{23} \dots  )
( \beta  _{31}  \circ \beta _{32} \circ \beta _{33} \dots  ) \
\dots,
 \end{equation*}
thus getting a proof for 
\cite[Condition C) on p.\ 281]{T}.
See the appendix.
In this respect, compare also \cite[Theorem 5]{CHL}. 
  The examples in \cite{dsh}
show that in some cases the bounds
given by Corollary \ref{coru} are  optimal or
close to be optimal.  

Applying 
Clause (4) above twice, we get a slightly different
proof, in comparison with \cite{CH},
that congruence modular varieties 
satisfy the \emph{Tolerance Intersection Property} (TIP)
$\Psi^*\Theta^* = (\Psi\Theta)^*$,
where $^*$ denotes transitive closure.
Again, Clause (4)
seems to give the best known bound 
for $\Psi \Theta^ \ell$. 
Notice that TIP has many important applications to congruence modular 
varieties.
See, e.~g., \cite{CHL}
for some examples and history.

\begin{proof}[Proof of Theorem \ref{=}]
If $\mathcal V$ is congruence distributive,
then $\mathcal V$ has indeed an alvin level $a(\mathcal V)$, by Theorem \ref{alvin}.
Since every congruence distributive variety is congruence modular, then
$\mathcal V$ has also a Gumm level $g(\mathcal V)$,
by Theorem \ref{gumm}.

A sequence of alvin terms is obviously also a sequence of  
Gumm terms, hence $g(\mathcal V) \leq a( \mathcal V)$.
On the other hand, suppose that $g(\mathcal V) = n $.
We will show that the inclusion 
\eqref{cd} holds, hence
 $a(\mathcal V) \leq n$,
by the equivalence of (1)  and (2) in Theorem \ref{alvin}.
Let $ (a, c) \in \alpha ( \beta \circ \gamma )$.
Since 
$\alpha ( \beta \circ \gamma ) \subseteq \alpha ( \beta + \gamma )$,
then,
 by congruence distributivity,
$ (a, c) \in \alpha  \beta + \alpha  \gamma $.
By Corollary \ref{coru}(1),
 $ (a, c) \in \alpha  \beta \circ \alpha \gamma 
\circ  {\stackrel{n}{\dots}} $
hence clause (2) in Theorem \ref{alvin}
holds. 
 \end{proof}  

Theorem \ref{=} can be proved by constructing a sequence of alvin 
terms $u_1, \dots u_n$ starting from a sequence
of Gumm terms   $t_1, \dots t_n$ and some other sequence of alvin terms of 
length not prescribed in advance.
 The procedure involves some deep nesting, naturally
leads to the introduction of terms of large a-rity and
might find further applications.
However, the simplest way to prove Theorem \ref{=}
seems the way we have presented here, using condition (2)
in Theorem \ref{alvin}.

A nontrivial $2$-Gumm term 
is a \emph{Maltsev term}, and it characterizes
  congruence permutable varieties.
A nontrivial $2$-alvin term is  a \emph{Pixley term}.   
Thus Theorem \ref{=} is a generalization of the
classical result that a congruence distributive variety $\mathcal V$ 
has a Pixley term if and only if  
$\mathcal V$ is
congruence permutable.

\begin{remark} \labbel{gummday}
An earlier characterization \cite{D} 
of congruence modularity
involves a sequence $u_0, \dots, u_r$ of quaternary \emph{Day terms}.
A variety with such a sequence is said to be
\emph{$r$-modular}. 
We shall not need the explicit definition of Day terms here.
See \cite{D,FMK,CV,ntcm,dsh,T} for full details.
From \cite{D} and \cite{G1} we get that a variety 
$\mathcal V$ has a sequence of Day terms
if and only if $\mathcal V$ has a sequence 
 of Gumm terms.
However, the relationship about the lengths of the above sequences is far
from being clear \cite{LTT}. 

Let the \emph{Day-to-Gumm function} DG be defined 
by setting DG($r$) to be the smallest $n$
such that every variety with a sequence $u_0, \dots, u_r$ of Day terms 
has a sequence 
$t_0, \dots, t_n$ of Gumm terms.
Let DG$_{\text{dist}}$
be defined in the same way, but restricting  ourselves to 
congruence distributive varieties.
In principle, it is possible that 
DG and DG$_{\text{dist}}$ are different functions.
The \emph{Gumm-to-Day functions} GD
and GD$_{\text{dist}}$
are defined symmetrically.

In \cite{G,LTT} it is shown 
that an argument from \cite{D}
carries over with minimal modifications in order to show
that $ \text{GD($n$)} \leq 2n-2$,
for $n \geq 2$.
It is proved in \cite{dsh} that this bound 
is optimal for $n$ even,
also when we restrict
ourselves to congruence distributive varieties.
Namely, we have
GD($n$) =
GD$_{\text{dist}}$($n$) =
$2n-2$, for $n$ even, $n \geq 2$.  
Less is known about  
 DG($r$).
The arguments 
in \cite{LTT}
show that 
 $ \text{DG($r$)} \leq r^2-r +1$
and in some special cases the value 
can be slightly improved.
However, the general problem
of evaluating DG($r$) seems completely open.
Provisional results about this and related problems appear
in \cite{ntcm}.
Notice that, as we mentioned,
  some results from 
\cite{ntcm} are improved here,
 so that \cite{ntcm} should be updated. 

Except for the evaluation of
DG($r$),
 the relationships among the
J{\'o}nsson, alvin, Day and Gumm levels have
been almost completely settled in 
\cite{FV,dsh} and the present work. 
In any case, using the above definitions 
and a result from \cite{LTT},
we get the following corollary,
showing that the Day level 
of a congruence distributive variety $\mathcal V$  
affects the distributivity levels of $\mathcal V$. 
 \end{remark}

\begin{theorem} \labbel{diLTT}
If $\mathcal V$ 
is an $r$-modular   congruence distributive variety,
then $\mathcal V$ is $r^2-r +2$-distributive.
More generally, if
$n= \text{\rm DG$_{\text{\rm dist}}$($r$)}$,
then $\mathcal V$ is
$n$-alvin and  $n+1$-distributive.
 \end{theorem} 

\begin{proof}
By definition, an 
$r$-modular   congruence distributive variety $\mathcal V$ 
is $n$-Gumm, for $n= \text{\rm DG$_{\text{\rm dist}}$($r$)}$.
By Theorem \ref{=}, $\mathcal V$ is $n$-alvin, hence
$n{+}1$-\brfrt distributive.
We have proved the second statement.
The first statement is then immediate from
the mentioned result from \cite{LTT}.    
\end{proof}

\textbf{Notice}. The present manuscript has been submitted to arXiv.
You might possibly find  copies
of this manuscript in other sites, 
such as  institutional repositories,
the author's web page or (deprived of  the present notice)
scholarly publications.

In case  you find a copy of the present version
(containing the present notice) on a site which requires 
a registration or a subscription, 
it is likely that such a  copy has 
been posted without the author's permission, consent or approval.
In any case, if some site contains a copy 
of the present version, it is not necessarily the case that the author is affiliated,
 supports or endorses in any way either the site
or the organization to which the site belongs.

\section*{Appendix} \labbel{app}

In this appendix we present the proofs of some results
only stated in the main body of the manuscript.

\begin{remark} \labbel{2}   
As already mentioned, the terms constructed in the proof
of Theorem \ref{reduction} have not minimal complexity.
For example, 
given any sequence 
 $t_0, \dots, t_n$ of Gumm  terms, the term 
$t_1$ satisfies 
(T$_2$).
Indeed, if 
$ a \mathrel { \Theta} b \mathrel {\Theta} c$, then 
$a= t_1(a,b,b)  \mathrel { \Theta} t_1(a,a,c)$. 
 \end{remark}

We now present 
a strengthening of
Theorem \ref{reduction}.
Not only we can replace the tolerance $\Theta$ 
in Theorem \ref{reduction} by a reflexive and compatible relation $R$.
We are allowed to replace \emph{each} occurrence of $\Theta$ 
by \emph{either} $R$ or $R ^\smallsmile$,   
where $R ^\smallsmile $ denotes the converse of $R$.

\begin{theorem} \labbel{redstr}
Suppose that $n \geq 1$
and $\mathcal V$ is a  variety 
with a sequence 
$t_0, \dots, t_n$ of Gumm (alvin) terms.
Then, for every $m \geq 1$,  $\mathcal V$ has a sequence 
$s_0, \dots, s_n$ of Gumm (alvin) terms
such that $s_1$ satisfies the 
following additional property.

  \begin{enumerate}    \item []
  \begin{enumerate}    \item    
[(A$_m$)] For every $\mathbf A \in \mathcal V$,
every reflexive compatible relation  $R$ on $\mathbf A$
and binary relations $X_1, \dots X_m$ 
such that, for every $j=1, \dots, m$,
either $X_j = R$ or 
$X_j = R ^\smallsmile $,
the following holds. 

If $a,c \in A$ and $ (a, c) \in X_1 \circ X_2
\circ \dots \circ  X_m  $, 
then $ a \mathrel { R } s_1 (a,a,c) $. 
 \end{enumerate} 
\end{enumerate} 
\end{theorem} 

In particular, if 
(A$_m$)
holds, then, by taking 
  $X_j = R$ for all $j$'s, we have that  
 $a \mathrel { R^m} c $ implies 
$ a \mathrel { R } s_1 (a,a,c) $.
By taking 
  $X_j = R ^\smallsmile $ for all $j$'s, we have that
 $a \mathrel { (R ^\smallsmile )^m} c $ implies 
$ a \mathrel { R } s_1 (a,a,c) $.
In general, 
(A$_m$)  asserts that $ a \mathrel { R } s_1 (a,a,c) $
follows from the assumption that 
$a$ and  $c$ are related by a chain of factors 
of length $m$, when each factor
is either $R$ or $R ^\smallsmile $.

Notice also that the term $s_1$ depends only on $m$, it 
does not depend on $R$, 
nor on the choice of the $X_j$'s.
  In particular, using the fact that 
$R ^{\smallsmile \smallsmile } =R$ and
 considering the converse of $R$
in place of $R$ in (A$_m$)
we also get $ a \mathrel { R ^\smallsmile  } s_1 (a,a,c) $,
that is, 
$ s_1 (a,a,c) \mathrel { R   }  a$.

\begin{proof}
By induction on $m$. If $m=1$, then the original term $t_1$ works, since  
if $a \mathrel { R } c $, then 
$a= t_1(a,a,a) \mathrel { R} t_1(a,a,c) $.
 This is the argument already presented in the proof of  
Theorem \ref{reduction}.
If $a \mathrel { R ^\smallsmile  } c $, that is,  
 $c \mathrel { R } a $, then 
$a= t_1(a,c,c) \mathrel { R } t_1(a,a,c) $.

In order to carry over the induction
we shall need two claims.
Let (A$^-_m$)
be defined like (A$_m$), but restricted to the cases
when the last factor $X_m $ is $  R$.
We first claim that if 
 $s_1$ satisfies (A$_m$),
then $s_1^{*}$
satisfies
(A$^-_{m+1}$),
where  $s_1^*$ is defined in \eqref{*}
in the proof of Theorem \ref{reduction}.  
Again, this is essentially the same argument.
Choose 
$X_1, \dots, X_m$
in some arbitrary way, let 
$X_{m+1} = R$ 
and set
$ X=  X_1 \circ X_2
\circ \dots \circ  X_m  $.
 If
$(a,c) \in X_1 \circ 
 \dots \circ  X_{m+1} $, then 
 $a \mathrel { X } b \mathrel {R } c  $,
for some $b$.
Then  
\begin{equation*}
    a= s_1 (a, \underline{s_1(a,a,b)}, s_1(a,a, \underline{b})) \mathrel { R } 
s_1 (a, \underline{a}, s_1(a,a, \underline{c})) =  s_1^*(a,a,c).
 \end{equation*}
We have $s_1(a,a,b) \mathrel { R } a$ 
by (A$_m$) and the remark shortly after the statement
of the present theorem.

Our second claim is that 
if 
(A$^-_m$) is satisfied by some term $r_1$, then  
(A$_m$)  is satisfied by  $r_1^*$,
Indeed, 
choose 
$X_1, \dots, X_m$
in an arbitrary way and, as above,  set
$ X=  X_1 \circ X_2
\circ \dots \circ  X_m  $.
Suppose that  $a \mathrel {X } c$. If $X_m= R$, then 
 $ a \mathrel { R } r_1 (a,a,c) $ by (A$^-_m$), hence
\begin{equation*}
a= r_1(a,a, \underline{a}) \mathrel { R} 
r_1(a,a, \underline{  r_1(a,a,c)}) = r_1^*(a,a,c).
 \end{equation*}    
On the other hand, if 
$X_m= R ^\smallsmile $, then we get  
 $ a \mathrel { R ^\smallsmile } r_1 (a,a,c) $
by applying  (A$^-_m$) with $R ^\smallsmile $
in place of $R$. Equivalently, 
 $ r_1 (a,a,c) \mathrel { R  } a  $.
 Then
\begin{equation*}
a= r_1(a,\underline{  r_1(a,a,c)}, r_1(a,a,c) )
 \mathrel { R} 
r_1(a,\underline{a},  r_1(a,a,c)) = r_1^*(a,a,c)
 \end{equation*}    
and our second claim follows.

Now suppose that
$s_0, \dots s_n$
is a sequence of Gumm terms and
$s_1$ satisfies (A$_m$).
By the arguments in the proof of Theorem \ref{reduction},
both  $s_0^*, \dots s_n^*$
and $s_0^{**}, \dots s_n^{**}$
are sequences of Gumm terms.
By our first claim,
$s_1^{*}$
satisfies
(A$^-_{m+1}$).
By our second claim
with $m+1$ in place of $m$
and  $s_1^{*}$ in place of 
$r_1$, we get that 
$s_1^{**}$
satisfies
(A$_{m+1}$),
thus the induction is complete.
\end{proof}

A  less general version of the following corollary
has been stated as an open problem in \cite[Condition C) on p.\ 281]{T}.
On the other hand, a more general result
is proved in \cite[Theorem 5]{CHL}.
However, the present proof can be used
to provide an explicit bound for
$\prod _{i \leq h} ( \beta _{i1} \circ  \beta _{i2} \circ \dots \circ  \beta _{ig} )$,
with possible repetitions among the $ \beta _{ij}$.

\begin{corollary} \labbel{condC}
For $ (\beta _{ij}) _{i\leq h, j \leq g}  $ congruences in
an algebra belonging to a congruence modular 
variety, we have 
\begin{multline*} 
\prod _{i \leq h} ( \beta _{i1} +  \beta _{i2} + \dots +  \beta _{ig} )
=
\\
\prod _{i \leq h} ( \beta _{i1} \circ   \beta _{i2} \circ \dots \circ  \beta _{ig} )
\circ \sum _{f : \{1, \dots, h\} \to  \{ 1, \dots, g \} } \beta _{1f(1)}  \beta _{2f(2)} \dots \beta _{if(i)}
\end{multline*}   
where $\prod$ denotes intersection.
 \end{corollary}

 \begin{proof}
Let $\Theta_i = 
\overline{\beta _{i1} \cup   \beta _{i2} \cup \dots \cup  \beta _{ig}}$, 
where $ \overline{R} $ denotes the smallest 
compatible relation containing $R$.
In particular,
$\Theta_i ^*= \beta _{i1} +  \beta _{i2} + \dots +  \beta _{ig}$. 
By  iterating 
Corollary \ref{coru}(4), we get
$\prod _{i \leq h} ( \beta _{i1} +  \beta _{i2} + \dots +  \beta _{ig} ) =
\Theta_1^*\Theta_2^*\dots \Theta_h^*=
(\Theta_1\Theta_2\dots \Theta_h)^*$.
Letting 
  $\Psi = \Theta_1\Theta_2\dots \Theta_h$,
we have $\Psi \subseteq 
\prod _{i \leq h} ( \beta _{i1} \circ   \beta _{i2} \circ \dots \circ  \beta _{ig} )$.
If $(a,c)$ belongs to the left-hand side
of the displayed formula, then
we get 
$a \mathrel { \Psi} s_1(a,a,c) $
by Theorem \ref{reduction}. 
Moreover
\begin{equation*}
(s_2(a,a,c),c) \in 
 \sum _{f : \{1, \dots, h\} \to  \{ 1, \dots, g \} } \beta _{1f(1)}  \beta _{2f(2)} \dots \beta _{if(i)},
  \end{equation*}     arguing as in the proof of 
Corollary \ref{coru}(3). 
\end{proof}

We now prove a more general version of 
Theorem \ref{=}. For simplicity, we state the version
dealing with congruence identities.
We also get a ``bilateral'' version.  

\begin{theorem} \labbel{bilateralga}
Suppose that $r \geq 1$
and, for each  $i$ with $ 1 < i \leq r$,
either $A_i = \alpha ( \gamma \circ \beta )$, or 
$A_i = \alpha  \gamma \circ  \alpha \beta $.
If $\mathcal V$ is a congruence distributive
variety  and $\mathcal V$ satisfies
\begin{equation}\labbel{S}\tag{S}
\alpha ( \beta \circ  \gamma )
\subseteq
\alpha (\gamma  \circ  \beta)  \circ A_2 \circ A_3 \circ  \dots
\circ  A _{ r -1} \circ A_r,
   \end{equation}    
then $\mathcal V$ satisfies  
\begin{equation}\labbel{S1}\tag{S1}
\alpha ( \beta \circ  \gamma  )
\subseteq
\alpha \gamma  \circ  \alpha \beta  \circ A_2 \circ A_3 \circ \dots
\circ  A _{r -1} \circ  A_r.
   \end{equation}    
If in addition $r \geq 2$ and  $A_r= \alpha ( \gamma \circ  \beta  )$,
then $\mathcal V$ satisfies  
\begin{equation}\labbel{S+}\tag{S$^+$}
\alpha ( \beta \circ  \gamma  )
\subseteq
\alpha \gamma  \circ  \alpha \beta  \circ A_2 \circ A_3 \circ \dots
\circ  A _{r -1} \circ \alpha \gamma  \circ \alpha \beta .
   \end{equation}
 \end{theorem} 

\begin{proof}
Assuming \eqref{S}, the standard arguments used to 
prove the equivalence of (1) and (2)  in Theorem \ref{alvin}
show that we have terms $t_1, \dots, t_n$ ($n= 2r$) satisfying
\eqref{a2} - \eqref{a5}, as well as 
$x=t_h(x,y,x)$, for a set of indices depending 
on the actual forms of the $A_i$'s.    
The arguments in the proofs of Theorem \ref{reduction}
show that, for every $m$,  we can have $t_1$ satisfying (T$_m$).  
Then the proof of Corollary \ref{coru}
shows that \eqref{S1} holds. 

If $A_r= \alpha ( \gamma \circ  \beta  )$, then, taking converses,
 \eqref{S1} is equivalent to  
$\alpha ( \gamma \circ  \beta  )
\subseteq
\alpha (\beta   \circ   \gamma )
  \circ A _{r -1} ^\smallsmile  \circ \dots  
\circ A_3 ^\smallsmile  \circ A_2 ^\smallsmile  
\circ  \alpha \beta   \circ  \alpha \gamma  $.
By applying the first statement with $\beta$ and $\gamma$ 
exchanged, we get
$\alpha ( \gamma \circ  \beta  )
\subseteq
\alpha \beta   \circ \alpha   \gamma 
  \circ A _{r -1} ^\smallsmile \circ \dots  
\circ A_3 ^\smallsmile  \circ A_2 ^\smallsmile  
\circ  \alpha \beta   \circ  \alpha \gamma  $,
hence \eqref{S+} follows by taking converses
again. 
\end{proof}

\end{document}